\input amstex
\documentstyle{amsppt}
\magnification=\magstep1
\hsize=5in
\vsize=7.3in
\TagsOnRight  
\topmatter
\title On a Conjecture of Kac-Wakimoto
%Fixed point resolutions.
\endtitle
\author Feng Xu \endauthor

\address{Department of Mathematics, University of Oklahoma, 601 Elm Ave,
Room 423, Norman, OK 73019}
\endaddress
\email{xufeng\@ math.ou.edu}
\endemail
\abstract{   
We prove a conjecture about mininmal index of certain representations of
Coset Algebraic Conformal Field Theories 
under certain conditions
as formulated previously by us. As a by-product, the Kac-Wakimoto 
Conjecture (KWC) which is related to the asympotics of the coset
characters is true under the same conditions.  The same idea in the proof
also proves a recent conjecture related to subfactors from conformal 
inclusions.}
\endabstract
\thanks     
I'd like to thank Professor Y. Kawahigashi for sending me [KLM]. \ \ \ \    
1991 Mathematics Subject Classification. 46S99, 81R10. 
\endthanks

\endtopmatter
\document   
            
\heading \S1. Introduction \endheading
 Let us first recall some definitions
from [X4]. \par
Let $G$ be a simply connected  compact Lie group
and let $H\subset G$ be a connected subgroup. Let $\pi^i$ be an irreducible
representations of $LG$ with positive energy at level $k$
\footnotemark\footnotetext{ When G is the direct product of simple groups,
$k$ is a multi-index, i.e., $k=(k_1,...,k_n)$, where $k_i\in \Bbb {N}$
corresponding to the level of the $i$-th simple group. The level of $LH$
is determined by the Dynkin indices of $H\subset G$. To save some writing
we write the coset as $H\subset G_k$.} 
on Hilbert space $H^i$.
Suppose when restricting to $LH$, $H^i$ decomposes as:
$$          
H^i = \sum_\alpha H_{i,\alpha} \otimes H_\alpha, 
$$ and  $\pi_\alpha$ are irreducible representations of $LH$ on
Hilbert space $H_\alpha$.  The set of $(i,\alpha)$ which appears in
the above decompositions will be denoted by $exp$. \par
To illustrate the new ideas in this paper we will focus on the case
when both $G$ and $H$ are  simply connected semisimple compact Lie groups of
type $A$, i.e., $G$ and $H$ are of the 
forms $SU(N_1)\times SU(N_2)\times...\times SU(N_n)$.
The ideas of this paper can be applied to all compact semisimple and
simply connected Lie
groups and we plan
to consider them  in  separate publications. \par 

We          
shall use $\pi^1$ (resp. $\pi_1$) \footnotemark\footnotetext{ This
is slightly different from the notation $\pi^0$ (resp. $\pi_0$) in [X4]:
it seems    
to be more appropriate since these representations correspond to identity
sectors.} to
denote the vacuum
representation of
$LG$ (resp. $LH$).
Let ${\Cal A}_{G/H}$ 
be the irreducible conformal precosheaf of the coset $G/H$ 
as defined in Prop. 2.2 of [X4].
The decompositions above naturally give rise to a class of
covariant representations of ${\Cal A}$, denoted by $\pi_{i,\alpha}$ or simply
$(i,\alpha)$.  By Th. 2.3 of [X4], $\pi_{1,1}$ is the vacuum representation
of ${\Cal A}_{G/H}$. \par

 Let us  
denote by $S_{ij}$ (resp. $\dot{S_{\alpha \beta}}$) the $S$ matrices of
$LG$        
(resp. $LH$) at level $k$ (resp. certain level of $LH$ determined by
the inclusion $H\subset G_k$) as 
defined on P. 264  of [Kac].  Define  \footnotemark\footnotetext {Our
$(j,\beta)$ corresponds to $(M,\mu)$ on P.186 of [KW], and it follows from
the definition (cf. \S2.1) and Cor. 2.10 of [GL]  
that $\langle (j,\beta), (1,1) \rangle$ is then equal
to $mult_M(\mu,p)$ which appears in 2.5.4 of [KW] if
the sector $(j,\beta)$ has finite index, and 
that our formula (1) is identical to 2.5.4 of [KW].}
$$          
b(i,\alpha) = \sum_{(j,\beta) } {S_{ij}} \overline{\dot S_{\alpha \beta}}
\langle (j,\beta), (1,1) \rangle \tag 1
$$          
Note the above summation is effectively over those $(j,\beta)$ such that 
$(j,\beta) \in exp$. The definition of $\langle x, y\rangle$ for
any two sectors $x,y$ can be found 
in \S2.1 before formula (0). \par
The Kac-Wakimoto Conjecture (KWC)  states that 
if $(i,\alpha) \in exp$, then $b(i,\alpha) >0$ (cf. Conj. 2.5 of [KW]).
This conjecture is related to the asymptotics of the coset characters
(cf. Th. B of [KW]).
\par
In \S2.4 of [X4] an even stronger conjecture, Conjecture 2 (C2) is
formulated. This conjecture states that the square root of the minimal
index (cf. \S2.1) 
of sector $(i,\alpha)$, denoted by $d_{(i,\alpha)}$, is given 
by 
$$
d_{(i,\alpha)} = \frac{b(i,\alpha)}{b(1,1)}.
$$ C2 is stronger than KWC since $d_{(i,\alpha)}\geq 1, b(1,1)> 0$.
\par 
In [X4], C2 is proved for a class of cosets, but the proof is based on
the known results about the branching rules, which already implies KWC.
The main improvement in this paper is a proof of C2 under general
conditions, without knowing details about the  branching rules
(cf. Th. 3.4). 
The power of the new ideas in this paper 
can be seen from the examples listed in \S3.1
after Cor. 3.5, where we show that KWC is true for infinite series of
cosets which do not seem to have been obtained by other methods. 
We also give a simple proof of  a conjecture (Conj. 7.1) in [BE3]
\footnotemark\footnotetext{ After this paper appears as 
preprint math.RT/9904098 in the net, the author is informed by
Prof. D. Evans that Conj. 7.1 is also proved in a forthcoming paper
(cf. [BEK1]) by 
different methods.}.\par
Let us describe the content of this paper in more details. In \S2 we 
collect some results from [Reh] and [KLM] which will be used in the 
proof of Prop. 3.1 in \S3. In particular the 
notion of modular matrices from [Reh],  and the notion of 
$\mu$-index from [KLM] are
introduced. In Lemma 2.2 we calculate the $\mu$-index of the
coset. Lemma 2.3 is an application of Prop. 3.1 of [BEK2] and
Prop. 3.1 of [BE4] to our setting.  
In Prop. 2.4 , we show that the global index
of the coset is the same as its $\mu$-index using Lemma 2.3. This
implies the non-degeneracy of the modular matrices for the coset
under the conditions of Prop. 2.4 
by Cor. 32 and Th. 38 of [KLM].  We also give a 
second proof of this result by using relative braidings first used 
in [X1] and studied in details in [BE3], [BEK1]. The 
$n$-regularity (cf. definition in \S2.2) of our coset (Cor. 2.5) follows from
Lemma 2.3 and Cor. 7 of [KLM].    

\par
The proof of Prop. 3.1 
contains one of the new ideas of this
paper, which is  a novel way of calculating the summation on the
right-hand side of the equation in (1) of Prop. 3.1.  
Prop. 3.1 is then used to
prove Cor. 3.2, Cor. 3.3 , Th. 3.4  and Cor. 3.5. The main results
of \S3,  Th. 3.4  and Cor. 3.5 have already been described at the
beginning of this introduction.  
The same idea in the proof of Prop. 3.1 is also used
in \S3.2 to give a proof of Conj. 7.1 in [BE3] (Th. 3.7). \par
As we already noted before,
all the Lie groups considered in this paper 
will be simply connected semisimple compact Lie group of
type $A$, i.e., groups of the form 
$SU(N_1)\times SU(N_2)\times...\times SU(N_n)$ unless stated otherwise. \par  
\heading 2. Preliminaries \endheading
\subheading {2.1 Genus 0 and 1 modular matrices }
Let us first recall some definitions from [X2].
Let $M$ be a properly infinite factor
and  $\text{\rm End}(M)$ the semigroup of 
 unit preserving endomorphisms of $M$.  In this paper $M$ will always
be the unique hyperfinite $III_1$ factors. 
Let $\text{\rm Sect}(M)$ denote the quotient of $\text{\rm End}(M)$ modulo 
unitary equivalence in $M$. We  denote by $[\rho]$ the image of
$\rho \in \text{\rm End}(M)$ in  $\text{\rm Sect}(M)$.\par
 It follows from
\cite{L3} and \cite{L4} that $\text{\rm Sect}(M)$, with $M$ a properly
infinite  von Neumann algebra, is endowed
with a natural involution $\theta \rightarrow \overline \theta $  ;  
moreover,  $\text{\rm Sect}(M)$ is
 a semiring with identity denoted by $id$ or $1$ when no confusion
arises. \par
If given a normal
faithful conditional expectation
$\epsilon:   
M\rightarrow \rho(M)$,  we define a number $d_\epsilon$ (possibly
$\infty$) by:
$$           
d_\epsilon^{-2} :=\text{\rm Max} \{ \lambda \in [0, +\infty)|
\epsilon (m_+) \geq \lambda m_+, \forall m_+ \in M_+
\}$$ (cf. [PP]).\par
 We define   
$$           
d = \text{\rm Min}_\epsilon \{ d_\epsilon |  d_\epsilon < \infty \}.
$$   $d$ is called the statistical dimension of  $\rho$. It is clear
from the definition that  the statistical dimension  of  $\rho$ depends
only         
on the unitary equivalence classes  of  $\rho$. 
The properties of the statistical dimension can be found in
[L1], [L3] and  [L4]. We will denote the statistical dimension of 
$\rho$ by $d_\rho$ in the following.  $d_\rho^2$ is called the minimal
index of $\rho$.\par
Recall from [X2] that  we denote 
by $\text{\rm Sect}_0(M)$ those elements
of           
$\text{\rm Sect}(M)$ with finite statistical dimensions.
For $\lambda $, $\mu \in \text{\rm Sect}_0(M)$, let
$\text{\rm Hom}(\lambda , \mu )$ denote the space of intertwiners from 
$\lambda $ to $\mu $, i.e. $a\in \text{\rm Hom}(\lambda , \mu )$ iff
$a \lambda (x) = \mu (x) a $ for any $x \in M$.
$\text{\rm Hom}(\lambda , \mu )$  is a finite dimensional vector 
space and we use $\langle  \lambda , \mu \rangle$ to denote
the dimension of this space.  $\langle  \lambda , \mu \rangle$
depends      
only on $[\lambda ]$ and $[\mu ]$. Moreover we have 
$$\langle \nu \lambda , \mu \rangle = 
\langle \lambda , \bar \nu \mu \rangle , 
\langle \nu \lambda , \mu \rangle 
= \langle \nu , \mu \bar \lambda \rangle 
$$ 
which follows from Frobenius 
duality (See \cite{L2} ).\par
Next we will recall some of the results of [Reh] (also cf. [FRS]) and
introduce
notations. \par
Let $\{[\rho_i], i\in I \}$ be a finite set  of
equivalence classes of irreducible
covariant representations of an irreducible conformal precosheaf
(cf. \S2.1 of [GL]) . 
Suppose this set is closed under
conjugation 
and composition. We will denote the conjugate of $[\rho_i]$ by
$[\rho_{\bar i}]$
and identity sector by $[1]$ if no confusion arises, and let 
$N_{ij}^k = \langle [\rho_i][\rho_j], [\rho_k]\rangle $. We will
denote by $\{T_e\}$ a basis of isometries in $\text {\rm
Hom}(\rho_k,\rho_i\rho_j)$. 
The univalence of $\rho_i$ (cf. P.12 of [GL]) will be denoted by
$\omega_{\rho_i}$. \par
Let $\phi_i$ be the unique minimal 
left inverse of $\rho_i$, define:
$$ 
Y_{ij}:= d_{\rho_i}  d_{\rho_j} \phi_j (\epsilon (\rho_j, \rho_i)^*
\epsilon (\rho_i, \rho_j)^*), \tag 0
$$ where $\epsilon (\rho_j, \rho_i)$ is the unitary braiding operator
 (cf. [GL] ). \par
We list two properties of $Y_{ij}$ (cf. (5.13), (5.14) of [Reh]) which
will be used in \S2.2:
$$
\align
Y_{ij} = Y_{ji} & = Y_{i\bar j}^* = Y_{\bar i \bar j} \tag 1 \\
Y_{ij} = \sum_k N_{ij}^k \frac{\omega_i\omega_j}{\omega_k} d_{\rho_k} \tag
2
\endalign
$$
 Define      
$\tilde \sigma := \sum_i d_{\rho_i}^2 \omega_{\rho_i}^{-1}$.
If the matrix $(Y_{ij})$ is invertible,
by Proposition on P.351 of [Reh] $\tilde \sigma$ satisfies
$|\tilde \sigma|^2 = \sum_i d_{\rho_i}^2$. 
Suppose $\tilde \sigma= |\tilde \sigma| \exp(i x), x\in {\Bbb R}$.
Define matrices 
$$           
S:= |\tilde \sigma|^{-1} Y, T:=  C Diag(\omega_{\rho_i})
\tag 3       
$$ where $C:= \exp(i \frac{x}{3}).$  Then these matrices satisfy the algebra:
$$           
\align       
SS^{\dag} & = TT^{\dag} =id, \tag 4  \\
TSTST&= S, \tag 5 \\
S^2 =\hat{C}, T\hat{C}=\hat{C}T=T, \tag 6
\endalign    
$$           
where $\hat{C}_{ij} = \delta_{i\bar j}$ is the conjugation matrix. Moreover
$$           
N_{ij}^k = \sum_m \frac{S_{im} S_{jm} S_{km}^*}{S_{1m}}. \tag 7
$$           
(7) is known as Verlinde formula. \par
We will refer the $S,T$ matrices
as defined in  (3)  as  {\bf genus 0 modular matrices} since
they are constructed from the fusions rules, monodromies and minimal
indices which can be thought as  genus 0 data associated to 
a Conformal Field Theory (cf. [MS]). \par
It follows from (7) and (4) that any irreducible representation
of the commutative ring generated by $i$'s is of the form
$i\rightarrow \frac{S_{ij}}{S_{1j}}$. \par
Now let us consider an example which verifies (1) to (7) above.
Let $G= SU(N)$. We denote $LG$ the group of smooth maps
$f: S^1 \mapsto G$ under pointwise multiplication. The
diffeomorphism group of the circle $\text{\rm Diff} S^1 $ is 
naturally a subgroup of $\text{\rm Aut}(LG)$ with the action given by 
reparametrization. In particular the group of rotations
$\text{\rm Rot}S^1 \simeq U(1)$ acts on $LG$. We will be interested 
in the projective unitary representation $\pi : LG \rightarrow U(H)$ that 
are both irreducible and have positive energy. This means that $\pi $ 
should extend to $LG\ltimes \text{\rm Rot}\ S^1$ so that
$H=\oplus _{n\geq 0} H(n)$, where the $H(n)$ are the eigenspace
for the action of $\text{\rm Rot}S^1$, i.e.,
$r_\theta \xi = \exp^{i n \theta}$ for $\theta \in H(n)$ and 
$\text{\rm dim}\ H(n) < \infty $ with $H(0) \neq 0$. It follows from 
\cite{PS} that for fixed level $K$ which
is a positive integer, there are only finite number of such 
irreducible representations indexed by the finite set
$$           
 P_{++}^{h}  
= \bigg \{ \lambda \in P \mid \lambda 
= \sum _{i=1, \cdots , N-1}
\lambda _i \Lambda _i , \lambda _i \geq 1\, ,
\sum _{i=1, \cdots , n-1}
\lambda _i < h \bigg \}
$$           
where $P$ is the weight lattice of $SU(N)$ and $\Lambda _i$ are the 
fundamental weights and $h=N+K$.  We will use 
$1$ to denote the trivial representation  of 
$SU(N)$. For $\lambda , \mu , \nu \in  P_{++}^{K}$, define
$$           
N_{\lambda \mu}^\nu  = \sum _{\delta \in P_{++}^{K} } \frac{S_{\lambda
\delta}      
S_{\mu \delta} S_{\nu \delta}^*}{S_{1\delta}} \tag 8
$$           
where $S_{\lambda\delta}$ is given 
by the Kac-Peterson formula:
$$           
S_{\lambda \delta} = c \sum _{w\in S_N} \varepsilon _w \exp
(iw(\delta) \cdot \lambda 2 \pi /n). \tag 9
$$           
Here  $\varepsilon _w = \text{\rm det}(w)$ and $c$ is a normalization 
constant fixed by the requirement that $(S_{\lambda\delta})$ 
is an orthonormal system. 
It is shown in \cite{Kac} P.288 that $N_{\lambda \mu}^\nu $ are
non-negative 
integers. Moreover, define $ Gr(C_K)$ 
to be the ring whose basis are elements 
of $ P_{++}^{K}$ with structure constants $N_{\lambda \mu}^\nu $.
  The natural involution $*$ on $ P_{++}^{K}$ is 
defined by $\lambda \mapsto \lambda ^* =$ the conjugate of $\lambda $ as 
representation of $SU(N)$.  All the irreducible representations of 
$Gr(C_K)$ are given by $\lambda \rightarrow
\frac{S_{\lambda\mu}}{S_{1\mu}}$ for some $\mu$. \par
The irreducible positive energy representations of $ L SU(N)$ at level
$K$ give rise to an irreducible conformal precosheaf ${\Cal A}_G$ 
and  its covariant representations by the results in  \S17 of [W2]. 
 ${\Cal A}_G$ is a collection of maps $I\in {\Cal I}\rightarrow 
{\Cal A}_G(I)$ from the proper intervals on a circle to von Neumann
algebras which satifsy conditions as defined in \S2 of [GL], and one can
also find the definitions of  covariant representations of  ${\Cal A}_G$
is \S2 of [GL].
We will use 
$\lambda$ to denote such   representations. \par
For $\lambda$ irreducible, the univalence $\omega_\lambda$ is given by
an explicit formula .  
Let us first 
define       
$$           
\Delta_\lambda = 
\frac {c_2(\lambda)}{K+N} \tag 10
$$ where $c_2(\lambda)$ is the value of
Casimir      
operator on representation of $SU(N)$ labeled by dominant weight
$\lambda$ (cf. 1.4.1 of [KW]).
 $\Delta_\lambda$ is usually called the conformal dimension.
Then         
we have      
$\omega_\lambda = \exp({2\pi i} \Delta_\lambda)$. 
Note that $\omega_\lambda=\omega_{\bar\lambda}$.  
\par

Define the central charge (cf. 1.4.2 of [KW]) 
$$           
C_G := \frac {K \text {\rm dim(G)}}{K+N} \tag 11
$$ and $T$ matrix as
$$           
T=diag(\dot\omega_\lambda) \tag 12
$$ where $\dot\omega_\lambda = \omega_\lambda exp (\frac{-2\pi i
C_G}{24}).    
$  By Th.13.8 of [Kac] $S$ matrix as defined in (9) and $T$ matrix
in (12) satisfy relation (4), (5) and (6). 
 Since $S,T$ matrix defined in (8) and
(11) are related to the modular properties of characters which are
related to Genus 1 data of CFT (cf. [MS]), we shall call them
{\bf genus 1 modular matrices.}
\par
By Cor.1 in \S34 of [W2],  The fusion ring generated by all
$\lambda \in   P_{++}^{(K)}$
is isomorphic to $ Gr(C_K)$, with structure constants $N_{\lambda
\mu}^\nu$ as defined in (8).
 One may therefore ask what are the $Y$ matrix
(cf. (0)) in this case. By using (2) and the formula for 
 $N_{\lambda 
\mu}^\nu$, a simple calculation shows:
$$           
Y_{\lambda \mu} = \frac{S_{\lambda \mu}}{S_{1\mu}},
$$ and it follows that $Y_{\lambda \mu}$ is nondegenerate, and $S,T$
matrices     
as defined in (3) are indeed the same $S,T$ matrix defined in (8) and
(11),        
which is a surprising fact. This fact is  refered to as
genus 0 modular matrices coincide with genus 1 modular matrices.  
If the analogue of Cor.1 in \S34 of 
[W2] is established for other types of simple and simply connected Lie 
groups, then this fact is also
true for other types of groups by the same argument. \par
\subheading {2.2 Nondegeneracy of the coset}
Let $H\subset G_k$ be as in the introduction. 
We will use ${\Cal A}_G, {\Cal A}_H$ to denote the irreducible
conformal precosheaves associated with $G$ and $H$ respectively
(see paragraph before (10) in \S2.1). Denote by 
${\Cal A}_{G/H}$ the irreducible
conformal precosheaves associated with the coset $H\subset G_k$
as defined in Prop. 2.2 of [X4]. \par
In [X4], certain rationality
results (cf. Th. 4.2) are proved for a class of coset $H\subset G_k$. 
A stronger rationality condition, $\mu$-rational or absolute rational,
is defined in \S3 of [KLM], and we will  recall these definitions.\par
Let  ${\Cal A}$ be  an irreducible conformal precosheaf  on 
a circle $S^1$. Two proper intervals $I_1, I_2$\footnotemark\footnotetext
{As in [GL] by an interval of the circle we mean an open connected
proper interval of the circle. If $I$ is such an interval then
$I'$ will denote the interior of the complement of $I$ in the circle.}
of the 
circle are said to be disjoint if $\bar I_1\cap \bar I_2=\emptyset.$
Denote by ${\Cal E_2}$ the set of two disjoint
intervals.  ${\Cal A}$ is said to be {\it split} if ${\Cal A}(I_1)\vee
{\Cal A}(I_2)$ is naturally isomorphic to the tensor product of
von Neumann algebras ${\Cal A}(I_1)\otimes {\Cal A}(I_2)$ for all
 $\{I_1, I_2\}\in {\Cal E_2}$. \par 
Let $\{I_1, I_2\}\in {\Cal E_2}$, and let
$I_3, I_4$ be the disjoint intervals such that $I_3\cup I_4$ is the
interior of the complement of $ I_1\cup I_2$ in $S^1$.   If the index   
$$[A(I_3')\cap A(I_4'): {\Cal A}(I_1)\vee
{\Cal A}(I_2)]
$$ is independent of $\{I_1, I_2\}\in {\Cal E_2}$, then this index is
called the  {\it $\mu$-index} of ${\Cal A}$, denoted by $\mu_{\Cal A}$.
 ${\Cal A}$ is said to be absolutely rational, or $\mu$-rational for
short, if ${\Cal A}$ is strongly additive with finite  $\mu$-index.\par
Recall that ${\Cal A}$ is {\it strongly additive} if ${\Cal A}(I_1)\vee 
{\Cal A}(I_2)= {\Cal A}(I)$ where $I_1\cup I_2$ is obtained 
by removing an interior point from $I$.  ${\Cal A}$ is said to
be {\it $n$-regular} if 
$ {\Cal A}(S^1-\{p_1,...,p_n\}) = {\Cal A}(S^1)$
for any $p_1,...,p_n \in S^1$. Note that if 
 ${\Cal A}$ is strongly additive, then 
 ${\Cal A}$ is $n$-regular for any $n$. \par
Here we  will make corrections on a statement concerning strong
additivity in [X4]. 
When $G$ is of type $A$, it 
is proved (cf. Remark on P. 504 of [W2] and [W4]) that
 ${\Cal A}_G$ is strongly additive. But as  
pointed out to us by Prof. A. Wassermann, the proof of this fact 
, Th. E in
P. 504  of [W2],  is not correct (This
does not affect the results of [W2] ).  
It follows that the proof of (1) of Lemma 2.1 and  the proof
of a remark on strong additivity at the end of \S2.1 in [X4]
which follows from the proof of Th. E in [W2] is not correct. 
We note that the correct proof of Th. E as remarked in [W2] applies
without change to give a proof of (1) of Lemma 2.1 in [X4]. The remark
on strong additivity at the end of \S2.1 in [X4] is used in \S4 of [X4]
only to ensure  the equivalence of local and global intertwinners
(also cf. \S2.3 of [BE1]), 
but under the condition of finite index and conformal
invariance,
the equivalence of local and global
intertwinners has been proved in Th. 2.3 of [GL]. Hence all the lemmas,
corollaries and theorems of [X2] hold without the strong additivity
assumption since one can instead use Th. 2.3 of [GL].  
However it will be interesting to prove strong additivity for
the net associated with the coset. See Cor. 2.5 for a positive
result in this direction. We will show in Prop. 2.4 how 
one can still get nondegeneracy of modular matrices for
the coset without knowing
strong additivity.\par
When ${\Cal A}$ has finite $\mu$-index, then ${\Cal A}$ has only
finitely many irreducible covariant representations with finite
index (cf. Th. 8 of [KLM] and the remark after it), denoted
by $\rho_i, i=1,...,n$. The global index of ${\Cal A}$, denoted by
$I_{{\Cal A}}$, is defined to be
$I_{{\Cal A}}:= \sum_{i} d_{\rho_i}^2$.  Note that by the proof of 
Th. 38 in [KLM] (also cf. \S3 of [X3]) 
$I_{{\Cal A}}\leq \mu_{{\Cal A}}$. 
\par
Denote by $\mu_G, \mu_H, \mu_{G/H}$ the $\mu$
indices of the irreducible conformal precosheaves
${\Cal A_G},  {\Cal A_H}$ and  ${\Cal A_{G/H}}$
associated with $H$, $G$ and the coset $H\subset G$ respectively.
The irreducible covariant representations of ${\Cal A_G},  {\Cal A_H}$ 
and  ${\Cal A_{G/H}}$ will be denoted by $i, \alpha$ and $x$ respectively.
 For simplicity the global index 
of ${\Cal A_G},  {\Cal A_H}$ 
and  ${\Cal A_{G/H}}$ will be denoted by $I_G, I_H$ and $I_{G/H}$
respectively. 
Note that as in \S1 $(1,1)$ will denote the vacuum representation of
 ${\Cal A_{G/H}}$. 
As in \S1, let $\pi^1$ denote the vacuum representation of
${\Cal A_G}$. 
Then one has a natural inclusion (cf. \S3 of [X4]) 
$$\pi^1({\Cal A_{G/H}}(I) \otimes {\Cal A_H}(I)) \subset
{\Cal A_G}(I)
$$  for any proper interval $I$ of a circle.  The coset
$H\subset G_k$ is {\it cofinite} if the above inclusion has finite index
(cf. \S3 of [X4]), and the square root of this index is denoted by 
$d(G/H)$. As noted in \S3 of [X4], $d(G/H)$ is independent of the choices
of $I$. By (3.1) of [X4], we have
$$
d(G/H)^2=\sum_{\alpha} d_{(1,\alpha)}d_\alpha \tag 13
$$ where $ d_{(1,\alpha)}$ and $d_\alpha$ are the statistical dimensions of
covariant representations $(1,\alpha)$ and $\alpha$ respectively.

We have:
\proclaim{Lemma 2.2}
Suppose that $H\subset G_k$ is cofinite, and
${\Cal A_H}$ and ${\Cal A_G}$ are  $\mu$-rational. 
Then ${\Cal A_{G/H}}$ is split and 
has finite $\mu$ index. In fact
$$
 \mu_{G/H} = \frac{d(G/H)^4 \mu_G}{\mu_H}         
$$
\endproclaim
\demo{Proof}
First note that  ${\Cal A_{G/H}}$ is split. This is a well known fact which
follows from the asymptotics of the growth of states in
the vaccum (cf. Th. B of [KW]) and [BAF]. For a simplified proof see
Prop. 2.3.1 of [X7] which follows from [W3]. 
Then the $\mu$-index of the tensor product of ${\Cal A_H}$ and 
${\Cal A_{G/H}}$  
is $\mu_H \mu_{G/H}$, and it follows 
from the proof of Prop.21 of [KLM] (also cf. the proof of Th. 3.5 in 
[X3]) that:
$$
\mu_H \mu_{G/H} = d(G/H)^4 \mu_G. 
$$  
It follows that
$$
 \mu_{G/H} = \frac{d(G/H)^4 \mu_G}{\mu_H}         
$$ is finite.
\enddemo
\hfill \qed
\par
Let us note the following interesting consequence of Lemma 2.2. For the 
diagonal inclusions of type A considered in \S2.2 of [X5], a direct
calculation using Lemma 2.2 shows that
$$
\mu_{G/H} = |\tilde \sigma|^2,
$$  where $\tilde \sigma$ is determined in (2) of Th. 2.3 in [X5].
By Th. 38 of [KLM], this shows that the irreducible sectors which are
determined in \S2.2 of [X5] are {\bf all} the irreducible sectors
of the coset theory, and by Cor. 7 of [KLM]
${\Cal A}_{G/H}$ is $n$-regular for
any $n$.  We will see  more general statements in Cor. 2.5 and Cor. 3.2.\par
The $\mu$-rationality of the irreducible confromal 
precosheaf associated with a type A group $G$
follows from the results of  [W2] and [X3]. In this case one has 
$I_G=\mu_G$.   
The proof of finite $\mu$ index and the calculation of
the index value in  [X3] are  based on the existence of a 
class of conformal inclusions which exist for all classical simply connected
Lie groups.\par   

%Here we present a different proof based on [X4] which
%may generalize to the case of exceptional groups.
%\proclaim{Proposition 2.3}
%If $\mu_G = \frac{1}{S_{11}^2}$ is true for $k=1$, then it is true for
%all $k\geq 1$.
%\endproclaim
%\demo{Proof}
%Let $H=SU(N)_k \subset G=SU(N)_{k-1} \times SU(N)_1$.
%Then 
%$$
%d(G/H)^2 = 
For the rest of this section, we assume $d(G/H)<\infty$. By Lemma 2.2 and 
the remark after Th. 8 of [KLM], ${\Cal A}_{G/H}$ has only finite number
of irreducible covariant representations. In fact the global index
$I_{G/H} \leq \mu_{G/H}$. Also note that the set of 
 irreducible covariant representations is closed under conjugation and 
composition (cf. [GL]). 
Note that if $x, \alpha$ are the  
irreducible covariant representations of  ${\Cal A}_{G/H}$ and
${\Cal A}_{H}$ respectively, then $x\otimes \alpha$ 
is an irreducible covariant representation
of ${\Cal A}_{G/H}\otimes{\Cal A}_{H}$. 
We can take  the finite set consisting of
$x\otimes \alpha$ where $x, \alpha$ are the  
irreducible covariant representations of  ${\Cal A}_{G/H}$ and
${\Cal A}_{H}$ respectively and define the $Y$-matrix as in (0) of
\S2.1. 
Then one has 
$$
Y_{x\otimes\alpha, y\otimes \beta} = Y_{xy} Y_{\alpha\beta},
$$ where $Y_{xy}, Y_{\alpha\beta}$ are the $Y$-matrix associated with
 the  set of 
irreducible covariant representations of  ${\Cal A}_{G/H}$ and
${\Cal A}_{H}$ respectively. $Y_{xy}$ will be referred to as the $Y$-matrix
of the coset. \par
In \S4.2 of [X4], the results of [X1] (also cf [BE1-2])
are applied to the net of
inclusions $\pi^1({\Cal A}_{G/H}(I)\otimes{\Cal A}_{H}(I))\subset
{\Cal A}_{G}(I).$  The key observation is that there is a ring
homorphism $x\otimes \alpha \rightarrow a_{x\otimes \alpha}$ with
certain remarkable properties first established in [X1].  
We will refer to \S4.2 of [X4] for the definition of 
$a_{x\otimes \alpha}$ and $\sigma_i$.  A useful property which follows from
Prop. 4.2  and (4) of Th. 4.1 of [X4] is
$$
\langle \sigma_i, a_{x\otimes \alpha}\rangle
= \langle (i,\alpha), x \rangle.
$$
\par
So the map $x\rightarrow a_{x\otimes 1}$ as defined
in [X4] is a ring isomorphism by   (1) of Prop.4.2 in [X4],
and  we have
$$
\langle x,y \rangle = \langle a_{x\otimes 1},a_{y\otimes 1} \rangle.
$$
We will use the notations of \S4.2 of [X4] and ideas of [X1]. 
We  denote the set of irreducible sectors of $ a_{x\otimes
\lambda}$  by $W$. \footnotemark\footnotetext
{Note this is slightly different from
the definition of vector space $V$ in \S3.1 of [X5],
and in fact $V\subset W$,
but we will
see in Cor. 3.2 that these two spaces coincide.}
Notice
$\sigma_i \in W$, and
 are referred to as ``special nodes" in \S3.4 of [X1].
The ring homomorphism $x\otimes \alpha\rightarrow a_{x\otimes \alpha}$
, up to a unitary equivalence, 
is called $\alpha$-induction in [BE1,2,3]. A dictionary between the
notations of [X1] and [BE1,2,3] can be found in \S2.1  of [X6].
If one choose the opposite braiding compared to the braiding
in the definition of $a_{x\otimes \alpha}$, one obtain 
$\tilde a_{x\otimes \alpha}$. Let $\tilde W$ be the set
of  irreducible subsectors of $\tilde a_{x\otimes \alpha}, \forall x, 
\alpha$. Let $\hat W$ be 
the set of  irreducible subsectors of $\tilde a_{x\otimes \alpha}
a_{y\otimes \beta}, \forall x, 
\alpha, y, \beta$. Let $W_0:= W\cap \tilde W$. For any finite set
$Z$ of irreducible sectors closed under
multiplication we denote by $I_Z:=\sum_{\theta\in Z}
d_\theta^2$, where $d_\theta$ is the statistical dimension of $\theta$
(cf. \S2.1). If no possible confusion arises, we will denote the 
vector space over ${\Bbb C}$ with basis $Z$ by $Z$. This vector
space is endowed with an inner product by extending the
bilinear form $\langle, \rangle$ on sectors (cf. \S2.1 )
linearily in the first variable and conjugate linearily in the
second variable. Note that
every sector of $Z$ gives rise to a linear operator acting on the
vector space $Z$ where the action is given by left multiplication.  

The following Lemma follows from Prop. 3.1 of
[BE4] and Prop. 3.1 of [BEK2]. We include the proof for our case.

\proclaim{Lemma 2.3}      
$$
I_W= I_{\tilde W}= \frac{I_H I_{G/H}}{d(G/H)^2},  I_{\hat W}=I_H I_{G/H},
I_{W_0}= \frac{I_H I_{G/H}}{d(G/H)^4}.
$$
\endproclaim
\demo{Proof}  
Each sector $a_{x\otimes \alpha}$ (and linear 
combinations of them) can be thought as an operator on
$W$ where the action is given by left multiplication. These operators
$a_{x\otimes \alpha}$ share a common eigenvector 
$d=\sum_{\lambda\in W} d_{\lambda} \lambda$ where $\lambda$ are
elements in the basis $W$, with eigenvalues 
$d_xd_\alpha$.  
Note that 
the matrix corresponding to $N=\sum_{x\otimes \alpha}a_{x\otimes \alpha}$  
on the basis $W$ is irreducible since each element of $W$ is a subsector
of  some $a_{x\otimes \alpha}$, and 
$d=\sum_{\lambda\in W} d_{\lambda} \lambda$ is also a 
Perron-Frobenius   eigenvector
of  $N$. Now define another vector 
$$
v:= \sum_{x,\alpha} d_x d_\alpha
a_{x\otimes \alpha}=\sum_{x,\alpha, \lambda} d_x d_\alpha
\langle a_{x\otimes \alpha}, \lambda\rangle \lambda.$$
Note that $v$ has positive entries under the basis $W$. 
Since
$$
\align
a_{y\otimes \beta} v= \sum_{x,\alpha} d_x d_\alpha
a_{y\otimes \beta}a_{x\otimes \alpha}
&= \sum_{x,z,\alpha,\delta} d_x d_\alpha N_{xy}^z N_{\alpha\beta}^{\delta}
a_{z\otimes \delta}\\
&= \sum_{z,\delta} d_y d_\beta d_\delta d_za_{z\otimes \delta}\\
&= d_y d_\beta  \sum_{z,\delta} d_z d_\delta a_{z\otimes \delta}
\endalign
$$ where we have used the homomorphism property of the map
$x\otimes \alpha\rightarrow a_{x\otimes \alpha}$. It follows that
$v$ is also an  Perron-Frobenius   eigenvector of $N$. So there exists
a positive constant $c$ such that
$v= c d$.  By
computing the statistical dimension we get
$I_H I_{G/H}= c I_W$. From $\langle v,id\rangle = c \langle d,id\rangle
$ where $id$ is the identity sector we get (using $d_{id}=1$)
$$
c= \sum_{x,\alpha} d_x d_{\alpha} \langle a_{x\otimes \alpha},id
\rangle =  \sum_{x,\alpha} d_x d_{\alpha} \langle x, (1,\alpha)
\rangle = \sum_{\alpha} d_{(1,\alpha)} d_{\alpha}
= d(G/H)^2
$$
where in the last equality we have used (13) in \S2. 
\par
Similarly if we define vectors $v^{-}= \sum_{x,\alpha} d_x d_\alpha
\tilde a_{x\otimes \alpha}, d^{-}= \sum_{\lambda\in \tilde W} d_\lambda$,
then a similar proof as above shows that
$v^{-}=c d^{-}$ and $I_H I_{G/H}= c I_{\tilde W}$.
This proves the first equation in the lemma. \par

To prove the second equation, consider $a_{x\otimes \alpha}
\tilde a_{y\otimes \beta}$ as operators on $\hat W$ with
the action given by multiplication on the left. Let
$\hat d= \sum_{\lambda\in \tilde W} d_\lambda$. $ \hat d$
is a commomn eigenvector of 
$a_{x\otimes \alpha}\tilde a_{y\otimes \beta}$ with eigenvalues
$d_x d_\alpha d_y d_\beta$. Define another vector
$$\hat v := 
\sum_{x,y,\alpha,\beta} d_x d_\alpha d_y d_\beta
a_{x\otimes \alpha}\tilde a_{y\otimes \beta}.
$$ One checks easily using the fact that 
as sectors $a_{x\otimes \alpha}\tilde a_{y\otimes \beta}=
\tilde a_{y\otimes \beta} a_{x\otimes \alpha}$ (cf. Lemma 3.3 of
[X1] or [BE2]) that $\hat v$ is also a commomn eigenvector of 
$a_{x\otimes \alpha}\tilde a_{y\otimes \beta}$ with eigenvalues
$d_x d_\alpha d_y d_\beta$. Let $\hat N:= 
\sum_{x,y,\alpha,\beta} a_{x\otimes \alpha}\tilde a_{y\otimes \beta}
.$ Then the matrix of $\hat N$ under the basis $\hat W$ is irreducible
since every irreducible sector in  $\hat W$ appears as a descendant of
some $a_{x\otimes \alpha}\tilde a_{y\otimes \beta}$. It follows that
$\hat d$, $\hat v$ are  Perron-Frobenius   eigenvectors of $\hat N$, and so
$ \hat v = \hat c \hat d$ for some positive constant $\hat c$. By
computing the statistical dimension we get:
$(I_H I_{G/H})^2 = \hat c I_{\hat W}$. From
$ \langle \hat v, id \rangle = \hat c \langle\hat d,id  \rangle$
we get (using $d_{id}=1$) 
$$
\align
\hat c=  \langle \hat v, id \rangle 
&= \sum_{x,y,\alpha,\beta} d_x d_\alpha d_y d_\beta
\langle a_{x\otimes \alpha}\tilde a_{y\otimes \beta},id \rangle \\
&=  \sum_{x,y,\alpha,\beta}  d_x d_\alpha d_y d_\beta
\langle a_{x\otimes \alpha}, \tilde a_{\bar y\otimes \bar\beta} \rangle 
\\
&= \sum_{x,y,\alpha,\beta}  d_x d_\alpha d_y d_\beta
\langle a_{x\otimes \alpha}, \tilde a_{ y\otimes \beta} \rangle 
\\  
%&= \sum_{x,y,\alpha,\beta}  Y_{x\otimes \alpha, (1,1)\otimes 1}
%\langle a_{x\otimes \alpha}, \tilde a_{ y\otimes \beta} \rangle 
%Y_{y\otimes \beta,(1,1)\otimes 1} 
%\\ 
&= \sum_{x,y,\alpha,\beta}  Y_{x\otimes \alpha, (1,1)\otimes 1}
\langle a_{x\otimes \alpha}, \tilde a_{ y\otimes \beta} \rangle
Y_{y\otimes \beta, (1,1)\otimes 1}\\
&=\sum_{x,y,\alpha,\beta} \langle id, a_{x\otimes \alpha}\rangle
Y_{x\otimes\alpha,y\otimes \beta} Y_{y\otimes\beta, (1,1)\otimes 1}
\endalign
$$ where in the last step we have used Th. 5.7 of [BEK1] as in the proof
of Prop. 3.1 in [BE4]. 
Since $ Y_{\alpha\beta} Y_{\beta 1}= I_H \delta_{\alpha,1}$,
$Y_{x\otimes\alpha,y\otimes \beta}= Y_{xy}Y_{\alpha\beta}$, 
we get 
$$
\align
\hat c= \sum_{x,y} \langle id, a_{x\otimes 1}\rangle
Y_{xy}Y_{y(1,1)} I_H
&= \sum_{x,y} \langle (1,1), x\rangle
Y_{xy}Y_{y1} I_H \\
&= \sum_{y} 
Y_{(1,1)y}Y_{y1} I_H = I_{G/H} I_H.
\endalign
$$ 
It follows that $ I_{\hat W}=I_H I_{G/H}$. \par
From the proof of the first equation above we have
$v=cd, v^{-}= c d^{-}$, and hence
$
\langle v, v^{-}\rangle = I_{W_0} I_W^2 / I_{W}^2.
$  But one can also compute directly that
$$
\langle v, v^{-}\rangle = \hat c = I_{G/H} I_H
$$
completing the proof of the third equality.
\enddemo
\qed
\par
Now we are ready to prove the following:
\proclaim{Proposition 2.4}
Suppose  $G$ and $H$ are simply connected semisimple compact
Lie groups of type A as noted in \S1.
Assume $H\subset G_k$ is also cofinite. 
Then:\par
(1) $I_{G/H}= \mu_{G/H}$;  \par
(2) The $Y$ matrix of the coset is nondegenerate.
\endproclaim
\demo{Proof}
Ad (1):
By the third equation in lemma 2.3
$I_{W_0} = \frac{I_H I_{G/H}}{d(G/H)^4}$. On the other hand since
$\sigma_i\in W_0, \forall i$ (cf. Lemma 3.5 of [X1] or [BE1]), it follows
by definition that
$I_G= \sum_i d_{\sigma_i}^2 \leq I_{W_0} = \frac{I_H I_{G/H}}{d(G/H)^4}$.
Since $I_G=\mu_G, I_H=\mu_H$, we get
$$
\mu_G \leq \frac{\mu_H I_{G/H}}{d(G/H)^4}.
$$  Note that
$I_{G/H}\leq \mu_{G/H}$, but by Lemma 2.2 
$$
\mu_G =\frac{\mu_H \mu_{G/H}}{d(G/H)^4}.
$$ It follows that all the $\leq$ above are $=$ and this proves (1). \par
Ad (2): This follows immediately from (1),  the second part of 
Th. 38 (note that strong additivity is not assumed) and Cor. 32 of 
[KLM]. \par
\enddemo
\qed \par
A different proof of (2) of Lemma 2.4 without using the results of
[KLM] can be given by using properties of relative braidings as follows.
Let $x$ be an irreducible
covariant representation of  ${\Cal A}_{G/H}$ which has trivial braidings
with every covariant representation of ${\Cal A}_{G/H}$, i.e., $\epsilon (x,y) 
\epsilon (y,x)= id, \forall y$. Then $x\otimes 1$ is a   
covariant representation of ${\Cal A}_{G/H}\otimes {\Cal A}_H$ which 
 has trivial braidings
with every covariant representation of ${\Cal A}_{G/H}\otimes {\Cal A}_H$. 
It follows by definition that $a_{x\otimes 1}= \tilde a_{x\otimes 1}$, and
so  $a_{x\otimes 1}\in W_0$. 
Since $I_{W_0}= I_G$, 
as sectors $ [a_{x\otimes 1}]= [\sigma_i]$ for some 
$i$ since $a_{x\otimes 1}$ is irreducible. 
Let $u$ be a unitary intertwinning operator such that  
$ \sigma_i = u a_{x\otimes 1}u^*$. For any $\sigma_j$, since
$ \sigma_j \prec a_{(j,\beta)\otimes \beta}$, we can choose an
isometry $u_1$ such that $\sigma_j= u_1^* a_{(j,\beta)\otimes \beta}
 u_1$. Note that by our assumption $x\otimes 1$ has trivial
braidings with $(j,\beta)\otimes \beta$, and by using $u,u_1$ and
the naturality 
of relative braidings (cf. Prop. 3.12, 3.15 of [BE3], also cf.
Lemma 2.2.3 of [X6]), we conclude that $ i$ and $j$
as irreducible covariant representations of ${\Cal A}_G$ have trivial 
braidings. This force $i$ to be identity, and it follows that
$ [a_{x\otimes 1}]= [id]$ which implies that $x$ is the identity
by Prop. 4.2 of [X4]. This proves (2) of Lemma 4.2 by the 
proposition in \S5 of [Reh].\par
\proclaim{Corollary 2.5}
Under the assumptions of Prop. 2.4 
${\Cal A}_{G/H}$ is $n$-regular for any $n\geq 1$.
\endproclaim
\demo{Proof}
This follows from Lemma 2.2, (1) of Lemma 2.4 and Cor. 7 of [KLM].
\enddemo
\qed \par 
\heading 3. {Kac-Wakimoto Conjecture} \endheading
\subheading {3.1. Conjecture 2 of [X4]} 
Let $H\subset G_k$ be as in \S1.  Throughout this section, we assume that  
$H\subset G_k$ is cofinite. \par
We will  denote by $ S,  T$
(resp. $\dot S, \dot T$) the  genus 0 modular matrices corresponding to
$G$ (resp. $H$). As we remarked at the end of \S2.1, they coincide
with the genus 1  modular matrices when $G,H$ are type $A$. \par
We also assume that  $H\subset G_k$ is not conformal, so the coset 
theory is non-trival (cf. Prop. 2.2 of [X4]).  
But see the remark after Prop. 3.1 for the
case of conformal inclusions.\par
By Prop. 2.4,  the
$Y$-matrix of the coset as defined in \S2.1 is non-degenerate, and 
we shall denote by
$\ddot S, \ddot T$ the corresponding genus 0 modular matrices. 
We will  denote by $ S,  T$
(resp. $\dot S, \dot T$) the  genus 0 modular matrices associated with
$G$ (resp. $H$). \par  
Throughout this section we will use  genus 0 
modular matrices only unless noted otherwise. \par    
As in \S2.2, we  denote the set of irreducible sectors of $ a_{x\otimes 
\alpha}$  by $W$. 
Notice
$\sigma_i \in W$, and 
 are referred to as ``special nodes" in \S3.4 of [X1].  
Since (cf. [X4] or [BE1-2]) 
$a_{\bar x\otimes \bar\alpha} = \bar a_{x\otimes \alpha},
\sigma_j
a_{x\otimes \alpha}
= a_{x\otimes \alpha} \sigma_j,$  the matrix corresponding to
multiplications on $W$ by  $\sigma_i, a_{x\otimes 1} 
,$  and $a_{1\otimes \alpha}$ are
commuting normal matrices, so they can be simultaneously diagonalized.
Note that  all the irreducible representations of the ring  generated
by $\alpha's$ are given by (cf. \S2.1 after (7))
$$           
\alpha \rightarrow \frac{\dot S_{\alpha \beta}}{\dot S_{1\beta}},
$$ and similarly for the ring generated by $\sigma_i'$s and 
$a_{x\otimes 1}$'s, with $\dot S$ replaced by $S$ and $\ddot S$ 
respectively.            
Assume  $ \{\psi^{(k,\delta,z;s)}\}$ are normalized
orthogonal   
eigenvectors of  the matrix corresponding to
multiplications on $W$ by  $\sigma_i, a_{x\otimes 1},$ and  
$a_{1\otimes \alpha}$ with eigenvalues
 $\frac {S_{i k}}{S_{1 k}}$, 
$\frac{\ddot S_{x z}}{S_{1z}}$ and 
$\frac{\dot S_{\alpha \delta }}{\dot S_{1\delta}}$ respectively, 
$s$ is an
index indicating the multiplicity of  $k,\delta,z$, and we denote by  
$(Exp1)$ the  set of $k,\delta,z;s$'s which appears in the set 
$ \{\psi^{(k,\delta,z;s)}\}$. 
Recall if a representation is denoted by $1$, it will always be the
vacuum representation.  
\proclaim{Lemma A}
We have:\par
(1)  The eigenvector $\psi^{(1,1,1;s)}$ is unique with multiplicity
$s=1$, and is  given by $\sum_a d_a a$, up to a positive constant; 
moreover $\sum_a d_a^2 = \frac{1}{|\psi_1^{(1,1,1;1)}|^2}$; 
\par
(2)  
$$
\sum_{(k,\delta,z;s)\in (Exp1)} 
\frac{S_{\bar i k}}{S_{1 k}}
\frac{\dot S_{\alpha \delta}}{S_{1 \delta}}\frac{\ddot S_{xz}}{\ddot S_{1
z}} |\psi_1^{(k,\delta,z;s)}|^2 = \langle \sigma_i, a_{x\otimes  \alpha}
\rangle;
$$
(3)  If
$$
\langle \sigma_j, a_{y\otimes  \beta}
\rangle \neq 0,
$$ then $\omega_y = \omega_j \omega_\beta^{-1}$; \par
(4) 
$$
\sum_{j,\beta,y} S_{1j} \overline{\dot{S_{1\beta}}}
\overline{\ddot{S_{1y}}} \langle \sigma_j, a_{y\otimes  \beta} 
\rangle =1.
$$
\endproclaim
\demo{Proof}
Ad (1):  Let $G= \sum_{x,\alpha} a_{x\otimes \alpha}$, then every element
of $W$ appears as an irreducible subsector of $G$, and so we have
$$
G_{ab} := \langle Ga, b \rangle = \langle G, b\bar a \rangle >0
$$ since $W$ is a ring. It follows that $(G_{ab})$ is an irreducible
matrix and has up to positive constant a unique Perron-Frobenius
eigenvector.  But the  vector  $\sum_a d_a a$ is an eigenvector of
$(G_{ab})$ by the properties of statistical dimensions with maximal
eigenvalue, and so  up to positive constant the  Perron-Frobenius
eigenvector  of $(G_{ab})$ is  $\sum_a d_a a$. Note that
 $\psi^{(1,1,1;s)}$ is an eigenvector of $(G_{ab})$ with 
 maximal eigenvalue, we must have $s=1$ and there exists a positive number 
$p$ such that 
$\psi^{(1,1,1;1)}_a = p d_a, \forall a$. The last part of
(1) now follows from $d_1=1$ and $\psi^{(1,1,1;1)}$ is a unit vector.
Note that there is an analogue statement in (3) of Th. 3.9 of [X1]. 
\par
Ad(2): Note $ \langle \sigma_i, a_{x\otimes  \alpha}
\rangle =  \langle  \sigma_{\bar i} a_{x\otimes  \alpha},1 
\rangle$, and (2) follows from the definitions. \par
Ad (3):  By (4) of Th.4.1 of [X4] (also cf. (2) of Prop.4.2 of [X4]) we
have
$$
\langle \sigma_j, a_{y\otimes  \beta}
\rangle = \langle (j,\beta),  y
\rangle, 
$$
so if    
$$
\langle \sigma_j, a_{y\otimes  \beta}
\rangle \neq 0,
$$ then $y$ appears as an irreducible sector of $ (j,\beta) $.  
Note that  the action of universal covering group ${\Bbb G}$ of
$PSL(2,{\Bbb R})$ (cf. Prop.2.2 of [GL])
on the Hilbert space $H_{(j,\beta)}$ induces an action on the
representation space $H_y$ corresponding to sector $y$, but  
the action of $2\pi$ in  ${\Bbb G}$ 
on the Hilbert space $H_{(j,\beta)}$ is given by a constant 
$\omega_j \omega_\beta^{-1}$, and it follows that the univalence 
$\omega_y$ of $y$ is given by
$\omega_y = \omega_j \omega_\beta^{-1}$. \par
Ad (4): 
By local equivalence (cf. Th. B in \S17 of [W2]), 
the minimal index of the subfactor 
$\pi^j (L_IH)'' \vee (\pi^j (L_IH)'\cap \pi^j (L_IG)'') \subset
 \pi^j (L_IG)''$ is independent of $j$, where $\pi^j$ is the
representation corresponding to $j$.  It follows from the properties of
statistical dimensions (cf. [L6]) that
$$
\sum_\beta d_{(j,\beta)} d_\beta = d_j^2 \sum_\beta d_{(1,\beta)}
d_\beta
= d_j^2 d(G/H)^2,
$$ where $d(G/H)^2= \sum_\beta d_{(1,\beta)}d_\beta$ and $d_{(j,\beta)}$
is the statistical dimension of the coset sector $(j,\beta)$. So we have:
$$
\align
\sum_{j,\beta,y} S_{1j} \overline{\dot{S_{1\beta}}}
\overline{\ddot{S_{1y}}} \langle \sigma_j, a_{y\otimes  \beta},
\rangle & = S_{11} \dot S_{11} \ddot S_{11}
\sum_{j,\beta,y} d_j d_\beta d_y
\langle (j,\beta), y \rangle \\
& =  S_{11} \dot S_{11} \ddot S_{11} \sum_{j,\beta} d_j d_{(j,\beta)} 
  =  S_{11} \dot S_{11} \ddot S_{11} \sum_j d_j^2 d(G/H)^2 \\
&=  \frac {d(G/H)^2 \dot S_{11} \ddot S_{11}}{S_{11}}.
\endalign
$$
Note that by our assumption $\mu_G = \frac{1}{S^2_{11}}, \mu_H =
\frac{1}{\dot S^2_{11}}$, and  $\mu_{G/H} =I_{G/H}
=\frac{1}{\ddot S^2_{11}}$ by (1) of Prop. 2.4,  
so it follows from Lemma 2.2 that
$$
\frac {d(G/H)^2 \dot S_{11} \ddot S_{11}}{S_{11}} =1,
$$ and the proof is complete. 
\enddemo
\qed
\proclaim{Proposition 3.1}
(1)
$$
\langle \sigma_i, a_{x\otimes  \alpha} \rangle
= \sum_{j,\beta,y} S_{ij} \overline{\dot{S_{\alpha\beta}}} 
\overline{\ddot{S_{xy}}} \langle \sigma_j, a_{y\otimes  \beta}
\rangle;    
$$
(2)
$$
\sum_x\langle \sigma_i, a_{x\otimes  \alpha} \rangle \ddot S_{xz}
= \sum_{j,\beta} S_{ij} \overline{\dot{S_{\alpha\beta}}}
 \langle \sigma_j, a_{ z\otimes \beta}
\rangle.
$$
\endproclaim
\demo{Proof} 
(2) Obviously follows from (1) and the 
unitarity of $\ddot S$,  so we just need to prove (1).
Use (cf. (2) and (3) of \S2.1) we have
$$
\ddot{S_{xy}} = \sum_w \langle xy,w\rangle \frac{\omega_x \omega_y}
{\omega_w} \ddot{S_{1w}}, 
$$ so the right hand side of Prop.3.1 is:
$$
\sum_{j,\beta,y,w}  S_{ij} \overline{\dot{S_{\alpha\beta}}}
\langle \sigma_j, a_{y\otimes  \beta}
\rangle \langle xy,w\rangle \frac{\omega_x^{-1} \omega_y^{-1}}
{\omega_w^{-1}} \ddot{S_{1w}} \tag *
$$  Note by (3) of Lemma A that if
$$
\langle \sigma_j, a_{y\otimes  \beta}
\rangle \neq 0,
$$ then $\omega_y = \omega_j \omega_\beta^{-1}$, so
we can substitute $\omega_j \omega_\beta^{-1}$ for 
$\omega_y$ in the above expression.  Now  take
the complex conjugate of both sides of  (2) of Lemma A with
$(i,\alpha, x)$ replaced by $(j,\beta,y)$ we have:
$$
\langle \sigma_j, a_{y\otimes \beta}
\rangle = \sum_{k,\delta,z;s} \frac{S_{jk}}{S_{1k}}
\frac{\dot S_{\bar\beta \delta}}{\dot S_{1\delta}}
\frac{\ddot S_{\bar y z}}{\ddot S_{1z}} 
|\psi_1^{(k,\delta,z,;s)}|^2,
$$ and we
shall plug this into
(*) and
we shall call this resulting
expression after the two substitutions above by (*) in the following.  
We first   sum over $j$ and $\beta$ in (*)  using (cf. (5) of \S2.1)
$$
\align
&\sum_j \omega_j^{-1} S_{ij} S_{jk} = C^3 \omega_i \omega_k S_{\bar i k},
\\
& \sum_\beta \omega_\beta \dot S_{\alpha \beta} \dot S_{\beta \delta} = \dot
C^{-3}
\omega_\alpha^{-1}
\omega_\delta ^{-1}  \dot S_{\alpha \delta}, 
\endalign
$$
and then   sum over $y$ in (*) using (cf. (2) of \S2.1)
$$
\sum_y \langle xy,w\rangle \frac{\ddot S_{\bar yz}}{\ddot S_{1z}}
= \frac{\ddot S_{xz}}{\ddot S_{1z}} \frac{\ddot S_{\bar wz}}{\ddot
S_{1z}},
$$ and finally sum 
over $w$ in (*) using (cf. (5) of \S2.1)
$$
\sum_w \ddot S_{1w} \omega_w \ddot S_{wz} =
\ddot C^{-3}
\omega_z^{-1} \ddot   
S_{1z}. 
$$ The right hand side of Prop. 3.1 is then 
$$
\sum_{(k,\delta,z;s)\in (Exp1)} \frac{C^3}{\dot C^3 \ddot C^3} \dot
\frac{\omega_i}{\omega_\alpha \omega_x} \dot
\frac{\omega_k}{\omega_\delta \omega_z}   
 \frac{S_{\bar i k}}{S_{1 k}}
\frac{\dot S_{\alpha \delta}}{S_{1 \delta}}\frac{\ddot S_{xz}}{\ddot S_{1
z}}
|\psi_1^{(k,\delta,z;s)}|^2.
$$ Set 
$$
i=1,\alpha=1,x=1,
$$ and use (4) of Lemma A, we have:
$$
\sum_{(k,\delta,z;s)\in (Exp1)} \frac{C^3}{\dot C^3 \ddot C^3} \dot
\frac{\omega_k}{\omega_\delta \omega_z}  
|\psi_1^{(k,\delta,z;s)}|^2 =1.
$$  By setting $i=1,\alpha=1,x=1$ in (2) of Lemma A we have
$$
\sum_{(k,\delta,z;s) \in (Exp1)}   |\psi_1^{(k,\delta,z;s)}|^2 =1.
$$ 
Since 
$$
|\frac{C^3}{\dot C^3 \ddot C^3}|=| 
\frac{\omega_k}{\omega_\delta \omega_z}| =1,
$$ 
and $|a+b| \leq |a| +|b|, \forall a, b\in {\Bbb C}$, we must have
$$
||\psi_1^{(1,1,1;1)}|^2 + \frac{\omega_k}{\omega_\delta \omega_z}
 |\psi_1^{(k,\delta,z;s)}|^2| = |\psi_1^{(1,1,1;1)}|^2 +
 |\psi_1^{(k,\delta,z;s)}|^2, \ \forall (k,\delta,z;s), 
$$

and since by (1) of Lemma A
$$
|\psi_1^{(1,1,1;1)}|^2 >0,
$$ we have
$$
\frac{\omega_k}{\omega_\delta \omega_z}
 |\psi_1^{(k,\delta,z;s)}|^2 = |\psi_1^{(k,\delta,z;s)}|^2,
$$
and it follows that  if
$$
|\psi_1^{(k,\delta,z;s)}| \neq 0
$$ for some $s$, then $\frac{\omega_k}{\omega_\delta \omega_z}=1$.

It follows that  
$ \frac{C^3}{\dot C^3 \ddot C^3} =1$ and  
we have proved that the RHS of Prop.3.1 is:
$$
\sum_{k,\delta,z;s}  \frac{S_{\bar i k}}{S_{1 k}}
\frac{\dot S_{\alpha \delta}}{S_{1 \delta}}\frac{\ddot S_{xz}}{\ddot S_{1
z}}
|\psi_1^{(k,\delta,z;s)}|^2,
$$ which is precisely the left hand side by (2) of Lemma A.
\enddemo
\qed
\par
Note that $ \frac{C^3}{\dot C^3 \ddot C^3} =1$ shows that
$\ddot C^3 = \frac{C^3}{\dot C^3} = \exp(-6\pi i \frac{(C_G -C_H)}{24})
$, where $C_G, C_H$ are central charges as defined in (11) of \S2.1.  This
matches with the result of [GKO] that the central charge of
the coset is  $C_G -C_H$.  For the special of diagonal cosets of type A,
this is  (2) of Th. 2.3 in [X5]
which is proved by different methods.\par  
Note that if $H\subset G_1$ is a conformal inclusion, then  
by a similar but simpler proof as above one can show the following:
$$
b_{i\alpha} = \sum_{j,\beta} S_{ij}\overline{\dot S_{\alpha \beta}} 
b_{j\beta},
$$ where $b_{i\alpha} \in {\Bbb N}$ are the branching coefficients.  This
is implied by (1) of Th. A on P. 185 of [KW]. \par
By setting $i=1,\alpha=1, x=1$ in (1) of Prop. 3.1 we have:
$$
\align
1& = \sum_{j,\beta,y} S_{1j} \overline{\dot S_{1\beta}}
\overline{\ddot S_{1y}} \langle \sigma_j, a_{y\otimes \bar \beta}
\rangle  \\
&= \sum_{s} \frac{1}{S_{11} \dot S_{11} \ddot S_{11}}
|\psi^{(1,1,1;s)}|^2 , 
\endalign
$$ so it follows from (1) of Lemma A that
$$
\sum_{a\in W} d_a^2 = \frac{1}{S_{11} \dot S_{11} \ddot S_{11}},
$$ where $a\in W$ means the summation over the basis of $W$ given
by irreducible sectors.\par

In [X5] we define $V$ to be the vector space whose basis are irreducible
subsectors of $\sigma_i a_{1\otimes \alpha}$, and $\psi^{(j,\beta,s)}$
are normalized eigenvectors of linear transformations on $V$ which
are multiplications by $\sigma_i $ and $a_{1\otimes \alpha}$. The 
proof of (1) of Lemma A applies in this case with $G=\sum_{i,\alpha}
\sigma_i a_{1\otimes \alpha}$ 
and we have
$$
\sum_{a\in V} d_a^2 = \frac{1}{S_{11} \dot S_{11} b^0(1,1)}, 
$$ where (cf. Prop. 3.1 of [X5]) 
$$
b^0(1,1):= \sum_{j,\beta} S_{1j} \dot S_{1\beta} \langle
\sigma_j,a_{1\otimes\beta} \rangle.
$$  Set $i=1,\alpha =1, z=1$ in (2) of Prop. 3.1
we have
$$
b^0(1,1) = \ddot S_{11}.
$$ \par
So we have proved that
$$
\sum_{a\in W} d_a^2 = \sum_{a\in V} d_a^2,
$$ and since $V\subset W$, we must have $V=W$. 
So for any irreducible sector $x$ of the coset, there exists 
$(j,\beta)$ such that $a_{x\otimes 1}$ is an irreducible subsector
of $\sigma_j a_{1\otimes \bar \beta}$.  Since
$$
\langle \sigma_j a_{1\otimes \bar \beta}, a_{x\otimes 1} \rangle
=\langle \sigma_j, a_{x\otimes \beta} \rangle
= \langle (j, \beta) ,  x \rangle, 
$$
we  have proved 
the following:
\proclaim{Corollary 3.2}
Every irreducible sector of the coset appears as an irreducible
subsector of some $(j,\beta) \in exp$.
\endproclaim
Note (1) of Th. 2.3 in [X5] follows from Cor. 3.2 above and Cor. 32 of
[KLM]. \par
Cor. 3.2 proves a stronger version of Conj. 1 in [X4] under the conditions
stated at the beginning of of this section.  It is interesting to note that
there is also a Vertex Operator Algebra (VOA) approach to the coset CFT in
[FZ]. In \S5 of [FZ] (also cf. P. 113 of [Kacv]) 
a coset VOA is defined and it is
conjectured that these coset VOA is rational (cf. [Z] for definitions).
For $(j,\beta) \in exp$, and $x$ a subsector of $ (j,\beta)$, 
let $H_{(j,\beta)}, H_x\subset H_{(j,\beta)} $ be the corresponding
Hilbert spaces of representations.  It is easy to see using \S2.3 of [X4]
that $H_x$ is also an irreducible representation of the coset VOA.  
This shows  Th. 4.3 of [X4] and Lemma 2.2 of [X5]
holds for the coset VOA in the case of diagonal coset of type A, which is
a result that has not been proved by using the theory of VOA so far.  
However  
to
use Cor. 3.2 to 
prove the rationality of the coset VOA, one has to show that 
any representation of this coset VOA admits a natural inner product so
similar  norm estimations as in  \S2.3 of [X4] can be carried out.
\footnotemark\footnotetext{ We'd like to thank Dr. Yongchang Zhu for a
discussion on this point.} \par 
 
More generally let us define 
$$
b^0(i,\alpha) := \sum_{j,\beta} S_{ij} \overline{\dot S_{\alpha\beta}}
\langle
\sigma_j,a_{1\otimes\beta} \rangle.
$$
Note as stated at the beginning of this section that all the 
$S,\dot S$ matrices above are  genus 0 $S, \dot S$ matrices.  So
$b^0(i,\alpha)$ is defined differently from $b(i,\alpha)$ in (1) of
\S1 where genus 1 $S, \dot S$ matrices are used.
\proclaim{Corollary 3.3}
The statistical dimension $d_{(i,\alpha)}$ of the sector $(i,\alpha)$ 
is given by
$$
d_{(i,\alpha)} = \frac{b^0(i,\alpha)}{b^0(1,1)}.
$$
\endproclaim
\demo{Proof}
By setting $i=1,\alpha=1, z=1$ in (2) of Prop. 3.1 we have
$b^0(1,1) = \ddot S_{11}$, and by setting $z=1$ in (2) of Prop. 3.1
we have
$$
\align
\frac{b^0(i,\alpha)}{b^0(1,1)}
 & = \sum_x\langle \sigma_i, a_{x\otimes  \alpha}
\rangle
 \frac{\ddot S_{x1}}{\ddot S_{11}} \\
&= \sum_x\langle \sigma_i, a_{x\otimes  \alpha}
\rangle d_x \\
&= \sum_x \langle (i,\alpha) , x
\rangle d_x \\ 
&= d_{(i,\alpha)}
\endalign
$$ which completes the proof of the Corollary.
\enddemo
\hfill \qed
\par
Note that if the genus 1 $S$ (resp. $\dot S$) matrix corresponding 
to $G$ (resp. $H$) coincide with the
genus 0 $S$  (resp. $\dot S$) corresponding
to $G$ (resp. $H$), then $b^0(i,\alpha)$ coincides with 
$b(i,\alpha)$ defined in (1) of \S1.
By Cor. 3.3, we have proved the following theorem:
\proclaim{Theorem 3.4}
Suppose  $G$ and $H$ are simply connected semisimple compact
Lie groups of type A as noted in \S1.
Assume $H\subset G_k$ is also cofinite. \par
Then Conjecture 2 in [X4] is true, i.e.,
the statistical dimension $d_{(i,\alpha)}$ of the coset sector
$(i,\alpha)$ 
is given by
$$
d_{(i,\alpha)} = \frac{b(i,\alpha)}{b(1,1)}.
$$
\endproclaim

Since $d_{(i,\alpha)}\geq 1, b(1,1) >0$, an immediate corollary of Th. 3.4
is the following:
\proclaim{Corollary 3.5}
Under the same conditions of Theorem 3.4, the Kac-Wakimoto Conjecture
is true, i.e., if $(i,\alpha) \in exp$, then
$b(i,\alpha) >0$.
\endproclaim
Let us mention some examples which satisfy the assumptions of
Th. 3.4. Take $H_l\subset G_1$ to be conformal inclusion where $H,G$ are
simply connected type $A$ Lie groups. Here is a list of such pairs:
$$           
\align       
{SU}(N)_{N-2} & \subset \ {SU} \left( \frac{N(N-1)}{2}
\right), \ \ N \geq 4 ;  \\
{SU}(N)_{N+2} & \subset \ {SU} \left( \frac{N(N+1)}{2}
\right),  \\
{SU}(M)_N \times {SU}(N)_M & \subset \ {SU}(NM).
\endalign
$$
Consider the coset $H_{lk} \subset G_k$ with $k\geq 2$.
By [W2] and [X3], $H_{lk}, G_k$ are $\mu$-rational and the genus 0
S-matrices coincide with the genus 1 S-matrix, and by (2) of Cor. 3.1 of
[X4] the coset $H_{lk} \subset G_k$ is cofinite, so the conclusion of
Th.3.4 and Cor.3.5 is true in these examples. To the best of our
knowledge, this is already a new result since
the branching rules ( the
set $exp$) is not known in
general\footnotemark\footnotetext{ The branching rules in the case of 
 conformal inclusions listed here are the main results of 
[LL] and [ABI] and are  by no means trivial.} for these examples, 
and even with the explict formula for
$exp$, the calculation of $b(i,\alpha)$ seems to be nontrivial in general. 
 
\subheading{ 3.2 Conjecture 7.1 of [BE3]}
We shall use the original settings of [X1]. 
Let $H_k \subset G_1$ be a conformal inclusion with both $G$ and $H$
being semisimple compact Lie groups of type $A$, 
and $k$ the Dynkin index of the
inclusion (cf. 
[KW]). We  use $i$ (resp. $\lambda$) to denote the irreducible
projective positive energy representation of loop group $LG$ (resp.
$LH$) at level 1 (resp. k) (cf.[PS]).\par  
Denote by $b_{i \lambda}$ the
branching
coefficients, i.e., when restricting to $LG$, $i$ decomposes as
$\sum_{\lambda} b_{i \lambda} \lambda$.  Denote by $S_{ij}$ (resp. 
$S_{\lambda \mu}$) the genus 1 S-matrices of  $LG$ (resp.
$LH$) at level 1 (resp. k) (cf. [Kac]). Recall $a_\lambda$
as defined on Page 372 of [X1], then we have (cf. Page 9 of [X2])
$$
b_{i \lambda} = \langle a_\lambda, \sigma_i \rangle.
$$

Let us first prove a lemma in the setting of [X1] which is an
analogue of Lemma 3.10 of [BE3].  The basic idea is already implicit
in the proof of Lemma 3.2 in [X1]. 
The proof depends on \S3 of [X1] and we refer the reader to [X1]
for unexplained notations.
\proclaim{Lemma 3.6}
If $\langle a_\lambda, \tilde a_\mu \rangle \neq 0$, then
$\omega_\lambda= \omega_\mu$.
\endproclaim
\demo{Proof}
Let $u\neq 0$ be in $Hom(a_\lambda, \tilde a_\mu)$. Then
$$
\rho(u) \in Hom (\gamma \lambda, \gamma \mu) =
Hom (\rho a_\lambda \bar\rho , \rho a_\mu \bar\rho), 
$$ 
so $u\in Hom(a_\lambda \bar\rho,  a_\mu\bar \rho)$, and by 
(1) of Th. 3.3 in [X1], we have 
$u\in Hom(a_\lambda ,  a_\mu)$. So we get:
$$
a_\mu(m) u = \tilde a_\mu (m) u, \forall m\in M.
$$  Set $m=w$, apply $\rho$ to both sides, and 
use the equation on P. 373
of [X1] we obtain
$$
\gamma(\sigma) w \rho(u) = \gamma (\tilde \sigma) w \rho(u).
$$ Multiply on the left by $v^*$ and use the equation on P. 369 of
[X1] we get:
$$
v^* \gamma(\sigma) w \rho(u) = d_\rho^{-1} \sigma \rho(u);
$$
$$
v^* \gamma(\tilde \sigma) w \rho(u) = d_\rho^{-1} \tilde \sigma \rho(u); 
$$ and so
$$
\sigma \rho(u) = \tilde \sigma \rho(u).
$$  Note this equation is an analogue of Lemma 3.6 of [BE3].\par
Now multiply both sides on the right by $v$, we have
$$
\sigma \rho(u) v  = \tilde \sigma \rho(u) v,
$$
hence 
$$
\sigma^{-1} \tilde \sigma \rho(u) v = \rho(u) v.
$$
Note that $\rho(u) v \in Hom (\lambda, \gamma \mu)$, and apply the
monodromy equation (cf. P. 359 of [X1] and use the fact that
the univalence of $\gamma$ is $1$) we get:
$$
\omega_\lambda \omega_\mu^{-1} \rho(u) v = \rho(u) v.
$$
To finish the proof  we just have to   show that
$\rho(u) v$  is not zero. Note by (3) of Prop. 2.6 of [X1]  we have
$\rho(u^*)= \gamma (u_1) w$ for some $u_1\in M$, and so
$ \rho(u) v = w^* \gamma (u_1^*) v =  w^* v  \gamma (u_1^*) 
= d_\rho^{-1} \gamma (u_1^*)$, so if 
$\rho(u) v=0$, then $\rho(u)=0$, and so $u=0$ contradicting our
assumption $u\neq 0$.  
\enddemo
\hfill \qed
\par
Let $U$ be the vector space with a basis which consists of irreducible
components of $a_\lambda \tilde a_\mu \sigma_i, \forall \lambda,\mu,i$.
$a_\lambda ,
\tilde a_\mu, \sigma_i$ acts on $U$ by multiplication, and since they are
normal
commuting matrices by (2) of Cor. 3.5 and Lemma 3.3 of [X1], they can
be simultaneously diagonalized, and suppose
$\{\psi^{(j,\lambda_1,\mu_1;s)} \}$ are normalized orthogonal eigenvectors 
of $a_\lambda ,
\tilde a_\mu,$ and $\sigma_i$ with eigenvalues 
$\frac{\dot S_{\lambda \lambda_1}}{ \dot S_{1 \lambda_1}},
\frac{\dot S_{\mu \mu_1}}{ \dot S_{1 \mu_1}},
$ and $\frac{ S_{i j}}{  S_{1 j}}$ respectively,
where $s$ is an index indicating the multiplicity of 
$j,\lambda_1,\mu_1$.  
\proclaim{Theorem 3.7}
$$
\langle \tilde a_\mu, a_\lambda \rangle = \sum_i b_{i\mu} b_{i\lambda}. 
$$
\endproclaim
\demo{Proof}
Let us calculate 
$$
\sum_{\lambda,\mu}\dot S_{1\mu} \langle \tilde a_\mu, a_\lambda \rangle
\dot S_{1\lambda} 
$$ as in the proof of Prop.3.1. By lemma 3.6, we have 
$$
\align
\sum_{\lambda,\mu}\dot S_{1\mu} \langle \tilde a_\mu, a_\lambda \rangle
\dot S_{1\lambda} &=
\sum_{\lambda,\mu}
\dot S_{1\mu}  \frac{\omega_\lambda}{\omega_\mu} \langle
\tilde a_\mu,
a_\lambda \rangle \dot S_{1\lambda} \\
&=\sum_{\lambda,\mu,i, \lambda_1,\mu_1;s} \dot S_{1\mu}
\omega_\mu^{-1}\dot S_{\bar \mu \mu_1} 
 \dot S_{1\lambda } \omega_\lambda\dot S_{\lambda \lambda_1} 
\frac{1}{\dot S_{1\lambda_1} \dot S_{1\mu_1}}
|\psi_1^{(i,\lambda_1,\mu_1;s)}|^2 \\
&= \sum_{i,\lambda_1,\mu_1;s} \omega_{\lambda_1}^{-1} \omega_{\mu_1}
|\psi_1^{(i,\lambda_1,\mu_1;s)}|^2,
\endalign
$$ 
where we have used (5) of \S2.1 in the last $=$.
It follows that
$$ 
\sum_{\lambda,\mu}\dot S_{1\mu} \langle \tilde a_\mu, a_\lambda \rangle
\dot S_{1\lambda}
\leq
\sum_{i,\lambda_1,\mu_1;s}                                       
|\psi_1^{(i,\lambda_1,\mu_1;s)}|^2 =\langle 1,1\rangle =1.
$$ On the other hand
$$ 
\align
\sum_{\lambda,\mu}\dot S_{1\mu} \langle \tilde a_\mu, a_\lambda \rangle
\dot S_{1\lambda} & \geq \sum_{\lambda,\mu,i}\dot S_{1\mu} \langle
\tilde
a_\mu, \sigma_i \rangle \langle \sigma_i,  a_\lambda \rangle
\dot S_{1\lambda} \\
&= \sum_{\lambda,\mu,i}\dot S_{1\mu} 
b_{i\mu }   b_{i\lambda}
\dot S_{1\lambda}  \\
& = \sum_{i} S_{i1} S_{1i} = 1,
\endalign
$$ 
where in the second $=$ we have used (a) of Th. A on  P. 185 of [KW].
 So we must have
$$
\sum_{\lambda,\mu}\dot S_{1\mu} \langle \tilde a_\mu, a_\lambda \rangle
\dot S_{1\lambda} = \sum_{\lambda,\mu,i}\dot S_{1\mu} 
b_{i\mu }   b_{i\lambda}
\dot S_{1\lambda},
$$  and since 
$$
\langle \tilde a_\mu, a_\lambda \rangle \geq  \sum_i b_{i\mu}
b_{i\lambda},
$$ and $\dot S_{1\mu}>0, \dot S_{1\lambda}>0$, we must have
$$
\langle \tilde a_\mu, a_\lambda \rangle =  \sum_i b_{i\mu}
b_{i\lambda}.
$$ 
\enddemo
\hfill \qed
\par
Th. 3.6 proves Conj. 7.1 of [BE3] is true. This together with Prop. 5.1 of
[BE3] show that the invariants of  the dual Jones-Wassermann subfactos
associated
with conformal inclusions are determined by the ring structure
generated by irreducible sectors of $a_\lambda \tilde a_\mu$,
thus removing the mystery expressed in the foonote on Page 393 
of [X1], where one can also find the first example of such ring.
\heading References \endheading
\roster
\item"{[ABI]}" D. Altschuler, M. Bauer and C. Itzykson, {\it The
branching rules of conformal embeddings},  Comm.Math.Phys.,
 {\bf 132}, 349-364
(1990). 
\item"[BE1]" J. B\"{o}ckenhauer, D. E. Evans,
{\it Modular invariants, graphs and $\alpha$-induction for
nets of subfactors. I.},  
Comm.Math.Phys., {\bf 197}, 361-386, 1998.
\item"[BE2]" J. B\"{o}ckenhauer, D. E. Evans,
{\it Modular invariants, graphs and $\alpha$-induction for
nets of subfactors. II.},  
Comm.Math.Phys., {\bf 200}, 57-103, 1999.
\item"[BE3]" J. B\"{o}ckenhauer, D. E. Evans,
{\it Modular invariants, graphs and $\alpha$-induction for
nets of subfactors. III.},  
Comm.Math.Phys., {\bf 205}, 183-228, 1999.
\item"[BE3]" J. B\"{o}ckenhauer, D. E. Evans,
{\it Modular invariants, graphs and $\alpha$-induction for
nets of subfactors. III.},  
Comm.Math.Phys., {\bf 205}, 183-228, 1999. Also see hep-th/9812110.
\item"[BE4]" J. B\"{o}ckenhauer, D. E. Evans,
{\it Modular Invariants from Subfactors: Type I Coupling 
Matrices and Intermediate Subfactors},  
math.OA/9911239.
\item"[BEK1]" J. B\"{o}ckenhauer, D. E. Evans, Y. Kawahigashi,
{\it On $\alpha$-induction, chiral generators and modular invariants
for subfactors}, Comm.Math.Phys., {\bf 208}, 429-487, 1999. Also
see math.OA/9904109.
\item"[BEK2]" J. B\"{o}ckenhauer, D. E. Evans, Y. Kawahigashi,
{\it Chiral structure of modular invariants for subfactors}, 
Comm.Math.Phys., {\bf 210}, 733-784, 2000.
\item"{[FRS]}" K.Fredenhagen, K.-H.Rehren and B.Schroer
,\par
{\it Superselection sectors with braid group statistics and 
exchange algebras. II}, Rev. Math. Phys. Special issue (1992), 113-157.
\item"{[FZ]}" I. Frenkel and Y. Zhu, 
{\it Vertex operator algebras associated to representations of 
affine and Virasoro algebras}, Duke Math. Journal (1992), Vol. 66, No. 1
, 123-168
\item"{[GL]}"  D.Guido and R.Longo, {\it  The Conformal Spin and
Statistics Theorem},  \par
Comm.Math.Phys., {\bf 181}, 11-35 (1996) 
\item"{[GKO]}" P. Goddard and D. Olive, eds., { \it Kac-Moody and
Virasoro algebras,} Advanced Series in Math. Phys., Vol 3, World
Scientific 1988.
\item"{[J]}" V. Jones, {\it Fusion en al\'gebres de Von Neumann et groupes
de lacets (d'apr\'es A. Wassermann),} Seminarie Bourbaki, 800, 1-20,1995.
\item"{[KLM]}" Y. Kawahigashi, R. Longo and M. M\"{u}ger,
{\it Multi-interval Subfactors and Modularity of Representations in
Conformal Field theory}, Preprint 1999, see also math.OA/9903104.
\item"{[KW]}"  V. G. Kac and M. Wakimoto, {\it Modular and conformal
invariance constraints in representation theory of affine algebras},  
Advances in Math., {\bf 70}, 156-234 (1988).
\item"{[Kac]}"  V. G. Kac, {\it Infinite Dimensional Lie Algebras}, 3rd
Edition,     
Cambridge University Press, 1990.
\item"{[Kacv]}"  V. G. Kac, {\it Vertex algebras for beginners}, AMS,
1998.
\item"{[Ka]}"  Y. Kawahigashi, {\it Classification of paragroup actions 
on subfactors}, Publ. RIMS, Kyoto Univ., {\bf 31} 481-517 (1995).
\item"{[L1]}"  R. Longo, Proceedings of International Congress of
Mathematicians, 1281-1291 (1994).
\item"{[L2]}"  R. Longo, {\it Duality for Hopf algebras and for
subfactors}, 
I, Comm. Math. Phys., {\bf 159}, 133-150 (1994).
\item"{[L3]}"  R. Longo, {\it Index of subfactors and statistics of
quantum fields}, I, Comm. Math. Phys., {\bf 126}, 217-247 (1989.
\item"{[L4]}"  R. Longo, {\it Index of subfactors and statistics of
quantum fields}, II, Comm. Math. Phys., {\bf 130}, 285-309 (1990).
\item"{[L5]}"  R. Longo, {\it An analog of the Kac-Wakimoto formula and
black hole conditional entropy,} gr-qc 9605073, to appear in
Comm.Math.Phys.
\item"{[L6]}"  R. Longo, {\it Minimal index and braided subfactors,
} J.Funct.Analysis {\bf 109} (1992), 98-112.
%\item"{[LD]}"  S. Doplicher and R. Longo, { Standard and split inclusions
%of von Neumann algebras.} Invent. Math. 73, 493 (1984)
%} J.Funct.Analysis {\bf 109} (1992), 98-112.
\item"{[LL]}"  F. Levstein, J.I. Liberati,   
{\it Branching Rules for Conformal Embeddings},   \par
Comm. Math. Phys., {\bf 173}, 1-16 (1995).
\item"{[LR]}"  R. Longo and K.-H. Rehren, {\it Nets of subfactors},
Rev. Math. Phys., {\bf 7}, 567-597 (1995).  
\item"{[MS]}" G. Moore and N. Seiberg, {\it Taming the conformal zoo},
Lett. Phys. B  {\bf 220}, 422-430, (1989).
\item"{[PP]}" M.Pimsner and S.Popa,
{\it Entropy and index for subfactors}, \par
Ann. Sci.\'{E}c.Norm.Sup. {\bf 19},
57-106 (1986). 
\item"[PS]" A. Pressley and G. Segal, {\it Loop Groups,} O.U.P. 1986.
\item"[Reh]" Karl-Henning Rehren, {\it Braid group statistics and their
superselection rules} In : The algebraic theory of superselection 
sectors. World Scientific 1990 
\item"[Reh2]" Karl-Henning Rehren, {\it Chiral observables and modular
invariants}, \par
 hep-th/9903262. 
\item"{[W1]}"  A. Wassermann, Proceedings of International Congress of
Mathematicians, 966-979 (1994).
\item"{[W2]}"  A. Wassermann, {\it Operator algebras and Conformal
field theories III},  Invent. Math. Vol. 133, 467-539 (1998)
\item"{[W3]}"  A. Wassermann, {\it  Operator algebras and Conformal
field theories}, preliminary notes of 1992.
\item"{[W4]}"  A. Wassermann, with contributions by V. Jones,{\it
Lectures on operator algebras and conformal field theory}, 
Proceedings of Borel Seminar, Bern 1994, to appear 
\item"[We]" H.Wenzl, {\it Hecke algebras of type A and subfactors},
Invent. Math. 92 (1988), 345-383
\item"[X1]" F.Xu, {\it   New braided endomorphisms from conformal
inclusions, } \par
Comm.Math.Phys. 192 (1998) 349-403.
\item"[X2]" F.Xu, {\it Applications of Braided endomorphisms from
Conformal inclusions,} 
Inter. Math. Res. Notice., No.1, 5-23 (1998),  also see q-alg/9708013,
and Erratum, Inter. Math. Res. Notice., No.8, (1998) 
\item"[X3]" F.Xu, {\it Jones-Wassermann subfactors for 
Disconnected Intervals}, \par
 q-alg/9704003, to appear in Comm. Contemp. Math.
\item"[X4]" F.Xu, {\it Algebraic coset conformal field theories},
 math.OA/9810035, to appear in Comm. Math. Phys.
\item"[X5]" F.Xu, {\it Algebraic coset conformal field theories II}, \par
 math.OA/9903096, Publ. RIMS, vol.35 (1999), 795-824.
\item"[X6]" F.Xu, {\it 3-manifold invariants from cosets},
math.GT/9907077.
\item"[Z]" Y. Zhu, {\it Modular invariance of characters of vertex 
operator algebras}, Journal AMS, 9 (1996), 237-302.
\endroster   
\enddocument